\documentclass{gtart}

\def\ifplaintex{\expandafter\ifx\csname documentclass\endcsname\relax}

\ifplaintex 
\else
\fi

\expandafter\ifx\csname beginpicture\endcsname\relax
\expandafter\ifx\csname documentclass\endcsname\relax
\input pictex \else
\input prepictex \input pictex \input postpictex \fi\fi

\def\gt{{\mathsurround=0pt\it $\cal G\mskip-2mu$eometry \&\ 
$\cal T\!\!$opology}}        

\def\gtp{{\mathsurround=0pt\it $\cal G\mskip-2mu$eometry \&\ 
$\cal T\!\!$opology $\cal P\!$ublications}}  

\def\volumenumber#1{\def\thevolumenumber{#1}}
\def\papernumber#1{\def\thepapernumber{#1}}
\def\volumeyear#1{\def\thevolumeyear{#1}}

\def\pagenumbers#1#2{\def\startpage{#1}\def\finishpage{#2}}
\def\published#1{\def\publishdate{#1}}
\def\proposed#1{\def\theproposer{#1}}
\def\seconded#1{\def\theseconders{#1}}
\def\received#1{\def\receiveddate{#1}}

\def\accepted#1{\def\accepteddate{#1}}
\def\asciititle#1{\def\theasciititle{#1}}
\def\covertitle#1{\def\thecovertitle{#1}}

\long\def\asciiabstract#1{\long\def\theasciiabstract{#1}}
\def\asciikeywords#1{\def\theasciikeywords{#1}}

\def\shorttitle#1{\def\theshorttitle{#1}}

\let\\\par
\let\thevolumenumber\relax\let\thepapernumber\relax
\let\thevolumeyear\relax\let\thesamplenumber\relax\let\startpage\relax
\let\finishpage\relax\let\publishdate\relax\let\receiveddate\relax
\let\reviseddate\relax\let\accepteddate\relax\let\theasciititle\relax
\let\thecovertitle\relax\let\theasciiauthors\relax
\let\theasciiabstract\relax\let\theasciikeywords\relax
\let\theasciiemail\relax\let\theshortauthors\relax\let\theshorttitle\relax

\long\def\maketitlep{   

\count0=\startpage

\gt\hfill      
\beginpicture
\setcoordinatesystem units <0.33truein, 0.33truein> point at 2.2 0.9
\setplotsymbol ({$\cal G$})
\plotsymbolspacing=9truept
\circulararc 315 degrees from 0 1 center at 0 0
\setplotsymbol ({$\cal T$})
\circulararc 315 degrees from 1 -1 center at 1 0
\endpicture
\break
{\small\ifx\thesamplenumber\relax 
Volume \else Sample
\fi\thevolumenumber\ (\thevolumeyear)
\startpage--\finishpage\nl
Published: \publishdate}
\vglue 0.5truein plus 0.4fil minus 0.1truein

{\parskip=0pt\leftskip 0pt plus 1fil\def\\{\par\smallskip}{\ifplaintex\large
\else\Large\fi\bf\thetitle}\par\medskip}   

\vglue 0pt plus 0.1fil 

{\parskip=0pt\leftskip 0pt plus 1fil\def\\{\par}{\sc\theauthors}
\par\medskip}

\vglue 0pt plus 0.1fil 

{\small\parskip=0pt\let\newline\\
{\leftskip 0pt plus 1fil\def\\{\par}{\sl\theaddress}\par}
\expandafter\ifx\theemail\relax    
\relax\else\vglue 5pt plus 0.02fil minus 2pt\def\\{\stdspace{\rm 
and}\stdspace} 
\cl{Email:\stdspace\tt\theemail}\fi
\ifx\theurl\relax                  
\relax\else\vglue 5pt plus 0.02fil minus 2pt\def\\{\stdspace{\rm 
and}\stdspace}
\cl{URL:\stdspace\tt\theurl}\fi\par}

\vglue 7pt plus 0.3fil minus 3pt

{\bf Abstract}
\vglue 5pt plus 0.1fil minus 2pt

\theabstract

\vglue 7pt plus 0.3fil minus 3pt

{\bf AMS Classification numbers}\quad Primary:\quad \theprimaryclass

Secondary:\quad \thesecondaryclass

\vglue 5pt plus 0.3fil minus 2pt

{\bf Keywords:}\quad \thekeywords

\vglue 10pt plus 0.5fil minus 5pt

{\small  Proposed: \theproposer\hfill Received: \receiveddate\nl
Seconded: \theseconders\hfill 
\ifx\reviseddate\relax                         
Accepted: \accepteddate                        
\else
Revised: \reviseddate                          
\fi}
\eject
}       

\let\maketitlepage\maketitlep
\let\maketitle\maketitlepage

\font\phead=cmsl9 scaled 950
\font=cmsl9 scaled 1050
\font\pnum=cmbx10 scaled 913
\font\lnum=cmbx10 
\font\pfoot=cmsl9 scaled 950
\font=cmsl9 scaled 1050
\ifplaintex
\headline{\vbox to 0pt{\vskip -4.5mm\line{\small\phead\ifnum
\count0=\startpage ISSN 1364-0380 (on line)
1465-3060 (printed) \hfill {\pnum\folio}\else\ifodd\count0\def\\{ }%
\ifx\theshorttitle\relax\thetitle\else\theshorttitle\fi\hfill{\pnum\folio}
\else\def\\{ and }{\pnum\folio}\hfill\ifx\theshortauthors\relax\theauthors
\else\theshortauthors\fi\fi\fi}\vss}}
\footline{\vbox to 0pt{\vglue 0mm\line{\small\pfoot\ifnum\count0=\startpage
\copyright\ \gtp\hfill\else
\gt, Volume \thevolumenumber\ (\thevolumeyear)\hfill\fi}\vss
}}
\else
\makeatletter
\def\@oddhead{{\small\lhead\ifnum\count0=\startpage ISSN 1364-0380 (on line)
1465-3060 (printed) \hfill {\lnum\number\count0}\else\ifodd\count0
\def\\{ }\ifx\theshorttitle\relax \thetitle \else\theshorttitle\fi\hfill
{\lnum\number\count0}\else\def\\{ and }{\lnum\number\count0}
\hfill\ifx\theshortauthors\relax 
\theauthors\else\theshortauthors\fi\fi\fi}}\def\@evenhead{\@oddhead}
\def\@oddfoot{\small\lfoot\ifnum\count0=\startpage\copyright\ \gtp\hfill\else
\gt, Volume \thevolumenumber\ (\thevolumeyear)\hfill\fi}
\def\@evenfoot{\@oddfoot}
\makeatother
\fi

\newwrite\gtoutfile
\long\gdef\makeheadfile{  
{\def\\{, }\def\s{ }
\immediate\openout\gtoutfile head.xxx
\immediate\write\gtoutfile{To: math@arxiv.org}
\immediate\write\gtoutfile{Subject: put or rep NNNNN:pppp}
\immediate\write\gtoutfile{--text follows this line--}
\immediate\write\gtoutfile{Proxy-for: \ifx\theasciiauthors\relax
\theauthors\else\theasciiauthors\fi\s<\ifx\theasciiemail\relax\theemail\else\theasciiemail\fi>}
\immediate\write\gtoutfile{\noexpand\\}
\immediate\write\gtoutfile{Authors: \ifx\theasciiauthors\relax
\theauthors\else\theasciiauthors\fi}
\immediate\write\gtoutfile{Title: \ifx\theasciititle\relax
\thetitle\else\theasciititle\fi}
\immediate\write\gtoutfile{Subj-class: GT or SG or MG etc}
\immediate\write\gtoutfile{MSC-class: \theprimaryclass\ifx\thesecondaryclass\relax\else, \thesecondaryclass\fi}
\immediate\write\gtoutfile{Journal-ref: Geom. Topol. \thevolumenumber
(\thevolumeyear) \startpage-\finishpage}
\immediate\write\gtoutfile{Comments: Published in Geometry and Topology at}
\immediate\write\gtoutfile{\s\s http://www.maths.warwick.ac.uk/gt/GTVol\thevolumenumber/paper\thepapernumber.abs.html}
\immediate\write\gtoutfile{\noexpand\\}
\immediate\write\gtoutfile{}
\ifx\theasciiabstract\relax
\immediate\write\gtoutfile{\theabstract}\else
\immediate\write\gtoutfile{\theasciiabstract}\fi
\immediate\write\gtoutfile{}
\immediate\write\gtoutfile{\noexpand\\}
\immediate\write\gtoutfile{}
\immediate\write\gtoutfile{<uuencoded .tar.gz file here>}
\immediate\write\gtoutfile{}
\immediate\closeout\gtoutfile}}  

\def\maketitlepage{\maketitlep\makeheadfile}
\let\maketitle\maketitlepage

\def\ifplaintex{\expandafter\ifx\csname documentclass\endcsname\relax}

\ifplaintex 
\else
\fi

\expandafter\ifx\csname beginpicture\endcsname\relax
\expandafter\ifx\csname documentclass\endcsname\relax
\input pictex \else
\input prepictex \input pictex \input postpictex \fi\fi

\def\gt{{\mathsurround=0pt\it $\cal G\mskip-2mu$eometry \&\ 
$\cal T\!\!$opology}}        

\def\gtp{{\mathsurround=0pt\it $\cal G\mskip-2mu$eometry \&\ 
$\cal T\!\!$opology $\cal P\!$ublications}}  

\def\volumenumber#1{\def\thevolumenumber{#1}}
\def\papernumber#1{\def\thepapernumber{#1}}
\def\volumeyear#1{\def\thevolumeyear{#1}}

\def\pagenumbers#1#2{\def\startpage{#1}\def\finishpage{#2}}
\def\published#1{\def\publishdate{#1}}
\def\proposed#1{\def\theproposer{#1}}
\def\seconded#1{\def\theseconders{#1}}
\def\received#1{\def\receiveddate{#1}}

\def\accepted#1{\def\accepteddate{#1}}
\def\asciititle#1{\def\theasciititle{#1}}
\def\covertitle#1{\def\thecovertitle{#1}}

\long\def\asciiabstract#1{\long\def\theasciiabstract{#1}}
\def\asciikeywords#1{\def\theasciikeywords{#1}}

\def\shorttitle#1{\def\theshorttitle{#1}}

\let\\\par
\let\thevolumenumber\relax\let\thepapernumber\relax
\let\thevolumeyear\relax\let\thesamplenumber\relax\let\startpage\relax
\let\finishpage\relax\let\publishdate\relax\let\receiveddate\relax
\let\reviseddate\relax\let\accepteddate\relax\let\theasciititle\relax
\let\thecovertitle\relax\let\theasciiauthors\relax
\let\theasciiabstract\relax\let\theasciikeywords\relax
\let\theasciiemail\relax\let\theshortauthors\relax\let\theshorttitle\relax

\long\def\maketitlep{   

\count0=\startpage

\gt\hfill      
\beginpicture
\setcoordinatesystem units <0.33truein, 0.33truein> point at 2.2 0.9
\setplotsymbol ({$\cal G$})
\plotsymbolspacing=9truept
\circulararc 315 degrees from 0 1 center at 0 0
\setplotsymbol ({$\cal T$})
\circulararc 315 degrees from 1 -1 center at 1 0
\endpicture
\break
{\small\ifx\thesamplenumber\relax 
Volume \else Sample
\fi\thevolumenumber\ (\thevolumeyear)
\startpage--\finishpage\nl
Published: \publishdate}
\vglue 0.5truein plus 0.4fil minus 0.1truein

{\parskip=0pt\leftskip 0pt plus 1fil\def\\{\par\smallskip}{\ifplaintex\large
\else\Large\fi\bf\thetitle}\par\medskip}   

\vglue 0pt plus 0.1fil 

{\parskip=0pt\leftskip 0pt plus 1fil\def\\{\par}{\sc\theauthors}
\par\medskip}

\vglue 0pt plus 0.1fil 

{\small\parskip=0pt\let\newline\\
{\leftskip 0pt plus 1fil\def\\{\par}{\sl\theaddress}\par}
\expandafter\ifx\theemail\relax    
\relax\else\vglue 5pt plus 0.02fil minus 2pt\def\\{\stdspace{\rm 
and}\stdspace} 
\cl{Email:\stdspace\tt\theemail}\fi
\ifx\theurl\relax                  
\relax\else\vglue 5pt plus 0.02fil minus 2pt\def\\{\stdspace{\rm 
and}\stdspace}
\cl{URL:\stdspace\tt\theurl}\fi\par}

\vglue 7pt plus 0.3fil minus 3pt

{\bf Abstract}
\vglue 5pt plus 0.1fil minus 2pt

\theabstract

\vglue 7pt plus 0.3fil minus 3pt

{\bf AMS Classification numbers}\quad Primary:\quad \theprimaryclass

Secondary:\quad \thesecondaryclass

\vglue 5pt plus 0.3fil minus 2pt

{\bf Keywords:}\quad \thekeywords

\vglue 10pt plus 0.5fil minus 5pt

{\small  Proposed: \theproposer\hfill Received: \receiveddate\nl
Seconded: \theseconders\hfill 
\ifx\reviseddate\relax                         
Accepted: \accepteddate                        
\else
Revised: \reviseddate                          
\fi}
\eject
}       

\let\maketitlepage\maketitlep
\let\maketitle\maketitlepage

\font\phead=cmsl9 scaled 950
\font=cmsl9 scaled 1050
\font\pnum=cmbx10 scaled 913
\font\lnum=cmbx10 
\font\pfoot=cmsl9 scaled 950
\font=cmsl9 scaled 1050
\ifplaintex
\headline{\vbox to 0pt{\vskip -4.5mm\line{\small\phead\ifnum
\count0=\startpage ISSN 1364-0380 (on line)
1465-3060 (printed) \hfill {\pnum\folio}\else\ifodd\count0\def\\{ }%
\ifx\theshorttitle\relax\thetitle\else\theshorttitle\fi\hfill{\pnum\folio}
\else\def\\{ and }{\pnum\folio}\hfill\ifx\theshortauthors\relax\theauthors
\else\theshortauthors\fi\fi\fi}\vss}}
\footline{\vbox to 0pt{\vglue 0mm\line{\small\pfoot\ifnum\count0=\startpage
\copyright\ \gtp\hfill\else
\gt, Volume \thevolumenumber\ (\thevolumeyear)\hfill\fi}\vss
}}
\else
\makeatletter
\def\@oddhead{{\small\lhead\ifnum\count0=\startpage ISSN 1364-0380 (on line)
1465-3060 (printed) \hfill {\lnum\number\count0}\else\ifodd\count0
\def\\{ }\ifx\theshorttitle\relax \thetitle \else\theshorttitle\fi\hfill
{\lnum\number\count0}\else\def\\{ and }{\lnum\number\count0}
\hfill\ifx\theshortauthors\relax 
\theauthors\else\theshortauthors\fi\fi\fi}}\def\@evenhead{\@oddhead}
\def\@oddfoot{\small\lfoot\ifnum\count0=\startpage\copyright\ \gtp\hfill\else
\gt, Volume \thevolumenumber\ (\thevolumeyear)\hfill\fi}
\def\@evenfoot{\@oddfoot}
\makeatother
\fi

\newwrite\gtoutfile
\long\gdef\makeheadfile{  
{\def\\{, }\def\s{ }
\immediate\openout\gtoutfile head.xxx
\immediate\write\gtoutfile{To: math@arxiv.org}
\immediate\write\gtoutfile{Subject: put or rep NNNNN:pppp}
\immediate\write\gtoutfile{--text follows this line--}
\immediate\write\gtoutfile{Proxy-for: \ifx\theasciiauthors\relax
\theauthors\else\theasciiauthors\fi\s<\ifx\theasciiemail\relax\theemail\else\theasciiemail\fi>}
\immediate\write\gtoutfile{\noexpand\\}
\immediate\write\gtoutfile{Authors: \ifx\theasciiauthors\relax
\theauthors\else\theasciiauthors\fi}
\immediate\write\gtoutfile{Title: \ifx\theasciititle\relax
\thetitle\else\theasciititle\fi}
\immediate\write\gtoutfile{Subj-class: GT or SG or MG etc}
\immediate\write\gtoutfile{MSC-class: \theprimaryclass\ifx\thesecondaryclass\relax\else, \thesecondaryclass\fi}
\immediate\write\gtoutfile{Journal-ref: Geom. Topol. \thevolumenumber
(\thevolumeyear) \startpage-\finishpage}
\immediate\write\gtoutfile{Comments: Published in Geometry and Topology at}
\immediate\write\gtoutfile{\s\s http://www.maths.warwick.ac.uk/gt/GTVol\thevolumenumber/paper\thepapernumber.abs.html}
\immediate\write\gtoutfile{\noexpand\\}
\immediate\write\gtoutfile{}
\ifx\theasciiabstract\relax
\immediate\write\gtoutfile{\theabstract}\else
\immediate\write\gtoutfile{\theasciiabstract}\fi
\immediate\write\gtoutfile{}
\immediate\write\gtoutfile{\noexpand\\}
\immediate\write\gtoutfile{}
\immediate\write\gtoutfile{<uuencoded .tar.gz file here>}
\immediate\write\gtoutfile{}
\immediate\closeout\gtoutfile}}  

\def\maketitlepage{\maketitlep\makeheadfile}
\let\maketitle\maketitlepage

\def\ifplaintex{\expandafter\ifx\csname documentclass\endcsname\relax}

\ifplaintex 
\else
\fi

\expandafter\ifx\csname beginpicture\endcsname\relax
\expandafter\ifx\csname documentclass\endcsname\relax
\input pictex \else
\input prepictex \input pictex \input postpictex \fi\fi

\def\gt{{\mathsurround=0pt\it $\cal G\mskip-2mu$eometry \&\ 
$\cal T\!\!$opology}}        

\def\gtp{{\mathsurround=0pt\it $\cal G\mskip-2mu$eometry \&\ 
$\cal T\!\!$opology $\cal P\!$ublications}}  

\def\volumenumber#1{\def\thevolumenumber{#1}}
\def\papernumber#1{\def\thepapernumber{#1}}
\def\volumeyear#1{\def\thevolumeyear{#1}}

\def\pagenumbers#1#2{\def\startpage{#1}\def\finishpage{#2}}
\def\published#1{\def\publishdate{#1}}
\def\proposed#1{\def\theproposer{#1}}
\def\seconded#1{\def\theseconders{#1}}
\def\received#1{\def\receiveddate{#1}}

\def\accepted#1{\def\accepteddate{#1}}
\def\asciititle#1{\def\theasciititle{#1}}
\def\covertitle#1{\def\thecovertitle{#1}}

\long\def\asciiabstract#1{\long\def\theasciiabstract{#1}}
\def\asciikeywords#1{\def\theasciikeywords{#1}}

\def\shorttitle#1{\def\theshorttitle{#1}}

\let\\\par
\let\thevolumenumber\relax\let\thepapernumber\relax
\let\thevolumeyear\relax\let\thesamplenumber\relax\let\startpage\relax
\let\finishpage\relax\let\publishdate\relax\let\receiveddate\relax
\let\reviseddate\relax\let\accepteddate\relax\let\theasciititle\relax
\let\thecovertitle\relax\let\theasciiauthors\relax
\let\theasciiabstract\relax\let\theasciikeywords\relax
\let\theasciiemail\relax\let\theshortauthors\relax\let\theshorttitle\relax

\long\def\maketitlep{   

\count0=\startpage

\gt\hfill      
\beginpicture
\setcoordinatesystem units <0.33truein, 0.33truein> point at 2.2 0.9
\setplotsymbol ({$\cal G$})
\plotsymbolspacing=9truept
\circulararc 315 degrees from 0 1 center at 0 0
\setplotsymbol ({$\cal T$})
\circulararc 315 degrees from 1 -1 center at 1 0
\endpicture
\break
{\small\ifx\thesamplenumber\relax 
Volume \else Sample
\fi\thevolumenumber\ (\thevolumeyear)
\startpage--\finishpage\nl
Published: \publishdate}
\vglue 0.5truein plus 0.4fil minus 0.1truein

{\parskip=0pt\leftskip 0pt plus 1fil\def\\{\par\smallskip}{\ifplaintex\large
\else\Large\fi\bf\thetitle}\par\medskip}   

\vglue 0pt plus 0.1fil 

{\parskip=0pt\leftskip 0pt plus 1fil\def\\{\par}{\sc\theauthors}
\par\medskip}

\vglue 0pt plus 0.1fil 

{\small\parskip=0pt\let\newline\\
{\leftskip 0pt plus 1fil\def\\{\par}{\sl\theaddress}\par}
\expandafter\ifx\theemail\relax    
\relax\else\vglue 5pt plus 0.02fil minus 2pt\def\\{\stdspace{\rm 
and}\stdspace} 
\cl{Email:\stdspace\tt\theemail}\fi
\ifx\theurl\relax                  
\relax\else\vglue 5pt plus 0.02fil minus 2pt\def\\{\stdspace{\rm 
and}\stdspace}
\cl{URL:\stdspace\tt\theurl}\fi\par}

\vglue 7pt plus 0.3fil minus 3pt

{\bf Abstract}
\vglue 5pt plus 0.1fil minus 2pt

\theabstract

\vglue 7pt plus 0.3fil minus 3pt

{\bf AMS Classification numbers}\quad Primary:\quad \theprimaryclass

Secondary:\quad \thesecondaryclass

\vglue 5pt plus 0.3fil minus 2pt

{\bf Keywords:}\quad \thekeywords

\vglue 10pt plus 0.5fil minus 5pt

{\small  Proposed: \theproposer\hfill Received: \receiveddate\nl
Seconded: \theseconders\hfill 
\ifx\reviseddate\relax                         
Accepted: \accepteddate                        
\else
Revised: \reviseddate                          
\fi}
\eject
}       

\let\maketitlepage\maketitlep
\let\maketitle\maketitlepage

\font\phead=cmsl9 scaled 950
\font=cmsl9 scaled 1050
\font\pnum=cmbx10 scaled 913
\font\lnum=cmbx10 
\font\pfoot=cmsl9 scaled 950
\font=cmsl9 scaled 1050
\ifplaintex
\headline{\vbox to 0pt{\vskip -4.5mm\line{\small\phead\ifnum
\count0=\startpage ISSN 1364-0380 (on line)
1465-3060 (printed) \hfill {\pnum\folio}\else\ifodd\count0\def\\{ }%
\ifx\theshorttitle\relax\thetitle\else\theshorttitle\fi\hfill{\pnum\folio}
\else\def\\{ and }{\pnum\folio}\hfill\ifx\theshortauthors\relax\theauthors
\else\theshortauthors\fi\fi\fi}\vss}}
\footline{\vbox to 0pt{\vglue 0mm\line{\small\pfoot\ifnum\count0=\startpage
\copyright\ \gtp\hfill\else
\gt, Volume \thevolumenumber\ (\thevolumeyear)\hfill\fi}\vss
}}
\else
\makeatletter
\def\@oddhead{{\small\lhead\ifnum\count0=\startpage ISSN 1364-0380 (on line)
1465-3060 (printed) \hfill {\lnum\number\count0}\else\ifodd\count0
\def\\{ }\ifx\theshorttitle\relax \thetitle \else\theshorttitle\fi\hfill
{\lnum\number\count0}\else\def\\{ and }{\lnum\number\count0}
\hfill\ifx\theshortauthors\relax 
\theauthors\else\theshortauthors\fi\fi\fi}}\def\@evenhead{\@oddhead}
\def\@oddfoot{\small\lfoot\ifnum\count0=\startpage\copyright\ \gtp\hfill\else
\gt, Volume \thevolumenumber\ (\thevolumeyear)\hfill\fi}
\def\@evenfoot{\@oddfoot}
\makeatother
\fi

\newwrite\gtoutfile
\long\gdef\makeheadfile{  
{\def\\{, }\def\s{ }
\immediate\openout\gtoutfile head.xxx
\immediate\write\gtoutfile{To: math@arxiv.org}
\immediate\write\gtoutfile{Subject: put or rep NNNNN:pppp}
\immediate\write\gtoutfile{--text follows this line--}
\immediate\write\gtoutfile{Proxy-for: \ifx\theasciiauthors\relax
\theauthors\else\theasciiauthors\fi\s<\ifx\theasciiemail\relax\theemail\else\theasciiemail\fi>}
\immediate\write\gtoutfile{\noexpand\\}
\immediate\write\gtoutfile{Authors: \ifx\theasciiauthors\relax
\theauthors\else\theasciiauthors\fi}
\immediate\write\gtoutfile{Title: \ifx\theasciititle\relax
\thetitle\else\theasciititle\fi}
\immediate\write\gtoutfile{Subj-class: GT or SG or MG etc}
\immediate\write\gtoutfile{MSC-class: \theprimaryclass\ifx\thesecondaryclass\relax\else, \thesecondaryclass\fi}
\immediate\write\gtoutfile{Journal-ref: Geom. Topol. \thevolumenumber
(\thevolumeyear) \startpage-\finishpage}
\immediate\write\gtoutfile{Comments: Published in Geometry and Topology at}
\immediate\write\gtoutfile{\s\s http://www.maths.warwick.ac.uk/gt/GTVol\thevolumenumber/paper\thepapernumber.abs.html}
\immediate\write\gtoutfile{\noexpand\\}
\immediate\write\gtoutfile{}
\ifx\theasciiabstract\relax
\immediate\write\gtoutfile{\theabstract}\else
\immediate\write\gtoutfile{\theasciiabstract}\fi
\immediate\write\gtoutfile{}
\immediate\write\gtoutfile{\noexpand\\}
\immediate\write\gtoutfile{}
\immediate\write\gtoutfile{<uuencoded .tar.gz file here>}
\immediate\write\gtoutfile{}
\immediate\closeout\gtoutfile}}  

\def\maketitlepage{\maketitlep\makeheadfile}
\let\maketitle\maketitlepage

\def\ifplaintex{\expandafter\ifx\csname documentclass\endcsname\relax}

\ifplaintex 
\else
\fi

\expandafter\ifx\csname beginpicture\endcsname\relax
\expandafter\ifx\csname documentclass\endcsname\relax
\input pictex \else
\input prepictex \input pictex \input postpictex \fi\fi

\def\gt{{\mathsurround=0pt\it $\cal G\mskip-2mu$eometry \&\ 
$\cal T\!\!$opology}}        

\def\gtp{{\mathsurround=0pt\it $\cal G\mskip-2mu$eometry \&\ 
$\cal T\!\!$opology $\cal P\!$ublications}}  

\def\volumenumber#1{\def\thevolumenumber{#1}}
\def\papernumber#1{\def\thepapernumber{#1}}
\def\volumeyear#1{\def\thevolumeyear{#1}}

\def\pagenumbers#1#2{\def\startpage{#1}\def\finishpage{#2}}
\def\published#1{\def\publishdate{#1}}
\def\proposed#1{\def\theproposer{#1}}
\def\seconded#1{\def\theseconders{#1}}
\def\received#1{\def\receiveddate{#1}}

\def\accepted#1{\def\accepteddate{#1}}
\def\asciititle#1{\def\theasciititle{#1}}
\def\covertitle#1{\def\thecovertitle{#1}}

\long\def\asciiabstract#1{\long\def\theasciiabstract{#1}}
\def\asciikeywords#1{\def\theasciikeywords{#1}}

\def\shorttitle#1{\def\theshorttitle{#1}}

\let\\\par
\let\thevolumenumber\relax\let\thepapernumber\relax
\let\thevolumeyear\relax\let\thesamplenumber\relax\let\startpage\relax
\let\finishpage\relax\let\publishdate\relax\let\receiveddate\relax
\let\reviseddate\relax\let\accepteddate\relax\let\theasciititle\relax
\let\thecovertitle\relax\let\theasciiauthors\relax
\let\theasciiabstract\relax\let\theasciikeywords\relax
\let\theasciiemail\relax\let\theshortauthors\relax\let\theshorttitle\relax

\long\def\maketitlep{   

\count0=\startpage

\gt\hfill      
\beginpicture
\setcoordinatesystem units <0.33truein, 0.33truein> point at 2.2 0.9
\setplotsymbol ({$\cal G$})
\plotsymbolspacing=9truept
\circulararc 315 degrees from 0 1 center at 0 0
\setplotsymbol ({$\cal T$})
\circulararc 315 degrees from 1 -1 center at 1 0
\endpicture
\break
{\small\ifx\thesamplenumber\relax 
Volume \else Sample
\fi\thevolumenumber\ (\thevolumeyear)
\startpage--\finishpage\nl
Published: \publishdate}
\vglue 0.5truein plus 0.4fil minus 0.1truein

{\parskip=0pt\leftskip 0pt plus 1fil\def\\{\par\smallskip}{\ifplaintex\large
\else\Large\fi\bf\thetitle}\par\medskip}   

\vglue 0pt plus 0.1fil 

{\parskip=0pt\leftskip 0pt plus 1fil\def\\{\par}{\sc\theauthors}
\par\medskip}

\vglue 0pt plus 0.1fil 

{\small\parskip=0pt\let\newline\\
{\leftskip 0pt plus 1fil\def\\{\par}{\sl\theaddress}\par}
\expandafter\ifx\theemail\relax    
\relax\else\vglue 5pt plus 0.02fil minus 2pt\def\\{\stdspace{\rm 
and}\stdspace} 
\cl{Email:\stdspace\tt\theemail}\fi
\ifx\theurl\relax                  
\relax\else\vglue 5pt plus 0.02fil minus 2pt\def\\{\stdspace{\rm 
and}\stdspace}
\cl{URL:\stdspace\tt\theurl}\fi\par}

\vglue 7pt plus 0.3fil minus 3pt

{\bf Abstract}
\vglue 5pt plus 0.1fil minus 2pt

\theabstract

\vglue 7pt plus 0.3fil minus 3pt

{\bf AMS Classification numbers}\quad Primary:\quad \theprimaryclass

Secondary:\quad \thesecondaryclass

\vglue 5pt plus 0.3fil minus 2pt

{\bf Keywords:}\quad \thekeywords

\vglue 10pt plus 0.5fil minus 5pt

{\small  Proposed: \theproposer\hfill Received: \receiveddate\nl
Seconded: \theseconders\hfill 
\ifx\reviseddate\relax                         
Accepted: \accepteddate                        
\else
Revised: \reviseddate                          
\fi}
\eject
}       

\let\maketitlepage\maketitlep
\let\maketitle\maketitlepage

\font\phead=cmsl9 scaled 950
\font=cmsl9 scaled 1050
\font\pnum=cmbx10 scaled 913
\font\lnum=cmbx10 
\font\pfoot=cmsl9 scaled 950
\font=cmsl9 scaled 1050
\ifplaintex
\headline{\vbox to 0pt{\vskip -4.5mm\line{\small\phead\ifnum
\count0=\startpage ISSN 1364-0380 (on line)
1465-3060 (printed) \hfill {\pnum\folio}\else\ifodd\count0\def\\{ }%
\ifx\theshorttitle\relax\thetitle\else\theshorttitle\fi\hfill{\pnum\folio}
\else\def\\{ and }{\pnum\folio}\hfill\ifx\theshortauthors\relax\theauthors
\else\theshortauthors\fi\fi\fi}\vss}}
\footline{\vbox to 0pt{\vglue 0mm\line{\small\pfoot\ifnum\count0=\startpage
\copyright\ \gtp\hfill\else
\gt, Volume \thevolumenumber\ (\thevolumeyear)\hfill\fi}\vss
}}
\else
\makeatletter
\def\@oddhead{{\small\lhead\ifnum\count0=\startpage ISSN 1364-0380 (on line)
1465-3060 (printed) \hfill {\lnum\number\count0}\else\ifodd\count0
\def\\{ }\ifx\theshorttitle\relax \thetitle \else\theshorttitle\fi\hfill
{\lnum\number\count0}\else\def\\{ and }{\lnum\number\count0}
\hfill\ifx\theshortauthors\relax 
\theauthors\else\theshortauthors\fi\fi\fi}}\def\@evenhead{\@oddhead}
\def\@oddfoot{\small\lfoot\ifnum\count0=\startpage\copyright\ \gtp\hfill\else
\gt, Volume \thevolumenumber\ (\thevolumeyear)\hfill\fi}
\def\@evenfoot{\@oddfoot}
\makeatother
\fi

\newwrite\gtoutfile
\long\gdef\makeheadfile{  
{\def\\{, }\def\s{ }
\immediate\openout\gtoutfile head.xxx
\immediate\write\gtoutfile{To: math@arxiv.org}
\immediate\write\gtoutfile{Subject: put or rep NNNNN:pppp}
\immediate\write\gtoutfile{--text follows this line--}
\immediate\write\gtoutfile{Proxy-for: \ifx\theasciiauthors\relax
\theauthors\else\theasciiauthors\fi\s<\ifx\theasciiemail\relax\theemail\else\theasciiemail\fi>}
\immediate\write\gtoutfile{\noexpand\\}
\immediate\write\gtoutfile{Authors: \ifx\theasciiauthors\relax
\theauthors\else\theasciiauthors\fi}
\immediate\write\gtoutfile{Title: \ifx\theasciititle\relax
\thetitle\else\theasciititle\fi}
\immediate\write\gtoutfile{Subj-class: GT or SG or MG etc}
\immediate\write\gtoutfile{MSC-class: \theprimaryclass\ifx\thesecondaryclass\relax\else, \thesecondaryclass\fi}
\immediate\write\gtoutfile{Journal-ref: Geom. Topol. \thevolumenumber
(\thevolumeyear) \startpage-\finishpage}
\immediate\write\gtoutfile{Comments: Published in Geometry and Topology at}
\immediate\write\gtoutfile{\s\s http://www.maths.warwick.ac.uk/gt/GTVol\thevolumenumber/paper\thepapernumber.abs.html}
\immediate\write\gtoutfile{\noexpand\\}
\immediate\write\gtoutfile{}
\ifx\theasciiabstract\relax
\immediate\write\gtoutfile{\theabstract}\else
\immediate\write\gtoutfile{\theasciiabstract}\fi
\immediate\write\gtoutfile{}
\immediate\write\gtoutfile{\noexpand\\}
\immediate\write\gtoutfile{}
\immediate\write\gtoutfile{<uuencoded .tar.gz file here>}
\immediate\write\gtoutfile{}
\immediate\closeout\gtoutfile}}  

\def\maketitlepage{\maketitlep\makeheadfile}
\let\maketitle\maketitlepage

\volumenumber{4}\papernumber{13}\volumeyear{2000}
\pagenumbers{397}{405}
\proposed{Robion Kirby}
\seconded{Wolfgang Metzler, Cameron Gordon}
\received{27 June 2000}
\accepted{3 November 2000}
\published{10 November 2000}

\usepackage{epsf,amssymb}

\newtheorem{thm}{Theorem}
\newtheorem{lemma}[thm]{Lemma}

\newcommand{\refine}{
$$\begin{picture}(140,160)  \footnotesize
    \put(-100,10)       {\epsfbox{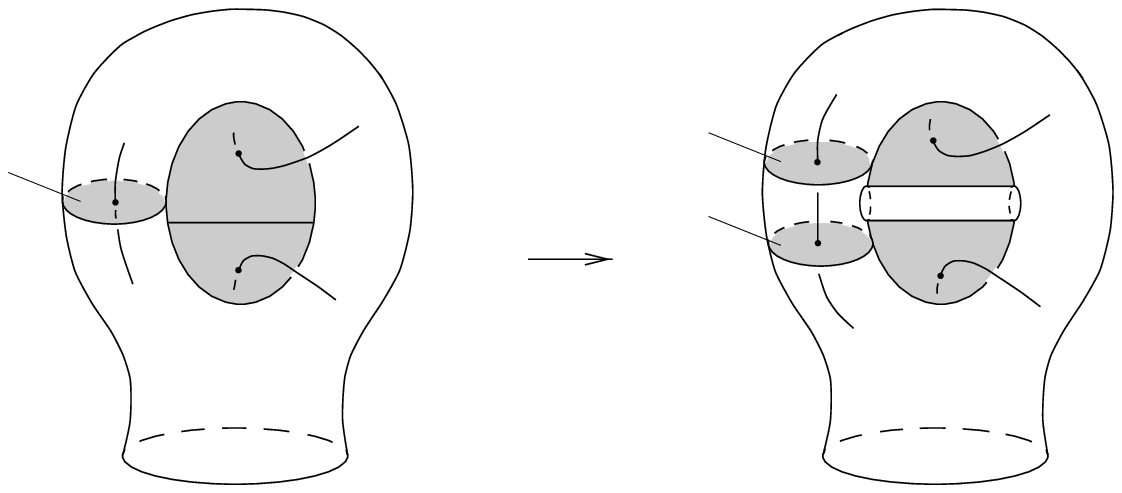}}
    \put(5,17)      {$\gamma$}
    \put(-90,50)    {$G^c$}
    \put(-70,115)   {${\Sigma}_i$}
    \put(-2,99)     {${\Sigma}_j$}
    \put(-5,70)     {${\Sigma}_k$}
    \put(-48,79)    {$\alpha$}
    \put(-111,100)  {$D$}
    \put(-48,92)    {$C$}
    \put(209,17)    {$\gamma$}
    \put(212,50)    {$G^c_{\rm split}$}
    \put(91,88)     {$D''$}
    \put(91,112)    {$D'$}
    \put(156,101)   {$C'$}
    \put(156,76)    {$C''$}
    \put(140,128)   {${\Sigma}_i$}
    \put(201,104)   {${\Sigma}_j$}
    \put(197,70)    {${\Sigma}_k$}
\end{picture}$$}

\newcommand{\exrefine}{
\unitlength 0.3mm
$$\begin{picture}(140,185)   \footnotesize
    \put(-120,10)       {\epsfxsize 4.45in \epsfbox{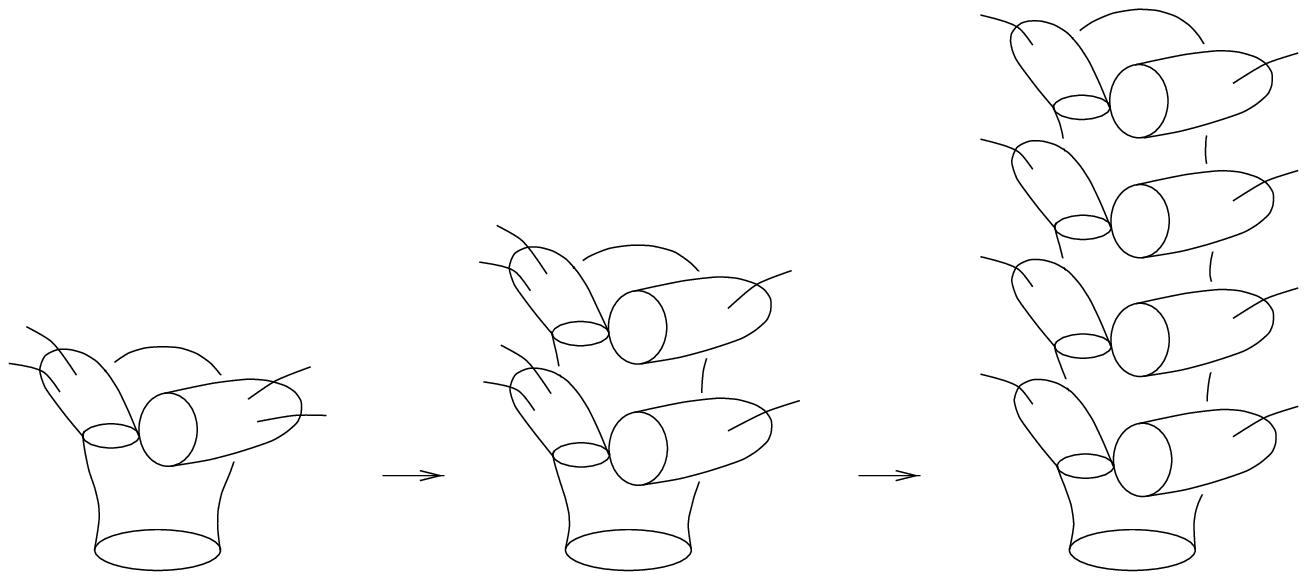}}
    \put(-130,83)    {${\Sigma}_1$}
    \put(-136,65)    {${\Sigma}_2$}
    \put(-27,71)     {${\Sigma}_3$}
    \put(-25,50)     {${\Sigma}_4$}
    \put(6,79)       {${\Sigma}_1$}
    \put(2,62)       {${\Sigma}_2$}
    \put(6,111)      {${\Sigma}_1$}
    \put(2,96)       {${\Sigma}_2$}
    \put(110,98)     {${\Sigma}_3$}
    \put(112,58)     {${\Sigma}_4$}
    \put(145,168)    {${\Sigma}_1$}
    \put(145,135)    {${\Sigma}_2$}
    \put(145,100)    {${\Sigma}_1$}
    \put(145,66)     {${\Sigma}_2$}
    \put(260,159)    {${\Sigma}_3$}
    \put(260,126)    {${\Sigma}_3$}
    \put(260,91)     {${\Sigma}_4$}
    \put(260,56)     {${\Sigma}_4$}
\end{picture}$$}

\newcommand{\grope}{
$$\begin{picture}(140,170)  \footnotesize
    \put(-35,0)       {\epsfbox{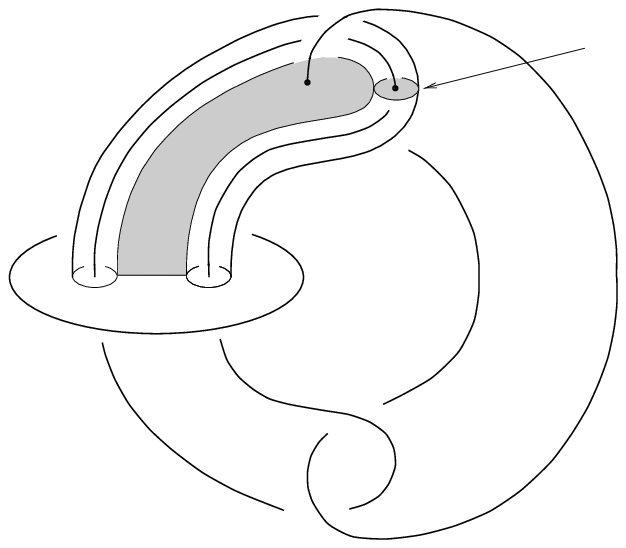}}
    \put(-15,125)    {$g^c$}
    \put(-50,75)     {$\gamma$}
    \put(-7,20)      {${\Sigma}_1$}
    \put(150,70)     {${\Sigma}_2$}
    \put(3,90)       {$C_2$}
    \put(136,142)    {$C_1$}
\end{picture}$$}

\begin{document}

\title{Exponential separation in ${\mathbf 4}$--manifolds}      
\shorttitle{Exponential separation in 4-manifolds}         
\asciititle{Exponential separation in 4-manifolds}         
\covertitle{Exponential separation in $4$--manifolds}      

\author{Vyacheslav S Krushkal}         
\newtheorem*{rem}{Remark}

\address{Department of Mathematics, Yale University\\New Haven, CT 06520-8283,
USA}
\email{krushkal@math.yale.edu}

\begin{abstract} 

We use a new geometric construction, grope splitting, to give a sharp
bound for separation of surfaces in $4$--manifolds. We also describe
applications of this technique in link-homotopy theory, and to the
problem of locating ${\pi}_1$--null surfaces in $4$--manifolds. In our
applications to link-homotopy, grope splitting serves as a geometric
substitute for the Milnor group.

\end{abstract}

\asciiabstract{We use a new geometric construction, grope splitting,
to give a sharp bound for separation of surfaces in 4-manifolds. We
also describe applications of this technique in link-homotopy theory,
and to the problem of locating pi_1-null surfaces in 4-manifolds. In
our applications to link-homotopy, grope splitting serves as a
geometric substitute for the Milnor group.}



\asciikeywords{4-manifolds, gropes, pi_1-null immersions, link homotopy}

\maketitlepage

Open problems in the classification theory of topological
four--manifolds, for ``large'' fundamental groups, have been
reformulated in terms of immersions of surfaces in $4$--manifolds, cf
\cite{Freedman0}, \cite{FQ}. Two related properties of immersed
surfaces are important in this discussion: disjointness, and vanishing
of the double point loops in the fundamental group of the ambient
$4$--manifold.  More precisely, these questions concern more general
$2$--complexes, {\em capped gropes}, that naturally arise in this
context. Capped gropes are assembled of several surface stages, capped
with disks with self-intersections. Each double point determines an
element of the fundamental group of the ambient $4$--manifold $M$, and
a central question is whether one can find a ${\pi}_1$--{\em null}
capped grope (so that each double point loop is contractible in $M$.)
A closely related question asks whether a collection of surfaces,
intersecting the caps of a given grope, may be pushed off it by a
homotopy, without creating intersections between different surfaces.
It follows from the work of Freedman--Teichner \cite{FT} that such
problems may be solved if the number of group elements (respectively,
the number of surfaces) is bounded by the exponential function $2^h-1$
in the grope height $h$.

In the present paper we describe a new construction, {\em grope
splitting}, which may be thought of as a tool for organizing
intersections between surfaces and capped gropes.  This construction
is used to give a new proof of the results on separation of surfaces,
and locating ${\pi}_1$--null capped gropes, mentioned above. The
argument in \cite{FT} relies on algebraic theory of link homotopy, and
is an indirect existence proof. Our proof is more transparent
geometrically, and it gives an explicit construction of the resulting
surfaces.  We also point out that the bound for separation of surfaces
is sharp.

This exponential result is one of the main ingredients of the theorem
\cite{FT} that the classification techniques ($4$--dimensional surgery
and $5$--dimensional $s$--cobor\-dism conjectures) hold in topological
category for fundamental groups of subexponential growth. A new
geometric proof of this theorem is presented in \cite{KQ}.  (Also see
the Appendix in that paper for a revised version of \cite{FT}.) The
conjectures for arbitrary fundamental groups remain open; the new
viewpoint presented here should be helpful in clarifying the problem.

In our applications to link-homotopy, the operation of grope splitting
replaces the Milnor group, used in the original proofs. Here grope
splitting is used to show that certain links are (colored)
link-homotopic. (Of course, Milnor group gives in general a more
precise algebraic information about links, but such generality is not
needed for the questions considered here.) In particular, we present a
new, geometric proof of the Grope Lemma, see Theorem \ref{grope lemma}
below.

We follow the terminology and notations of \cite{FT}. In particular,
$g$ denotes a grope (the underlying $2$--complex), while the capital
letter $G$ indicates the use of its untwisted $4$--dimensional
thickening.  The body of a capped grope $g^c$ is denoted by $g$. We
refer the reader to \cite{FQ} for definitions and a discussion of the
properties of gropes. The operations that are used extensively in this
paper (justifying the informal name for theorem \ref{ecp} below) are
{\em contraction}, sometimes also referred to as {\em symmetric
surgery}, and {\em pushoff}, which are described in detail in
\cite[section 2.3]{FQ}.  We remark that these operations are suited
perfectly for the purpose of separating surfaces in $4$--manifolds at
the expense of introducing self-intersections.

\begin{thm}[Exponential separation]\label{surfaces}
Let $(g^c,{\gamma})$ be a capped grope of height $h$, properly immersed in a 
$4$--manifold $M$, and let ${\Sigma}_1,\ldots,{\Sigma}_{2^h-1}$ 
be properly immersed surfaces in $M$ which are pairwise disjoint, and are also disjoint 
from the body of $g^c$. Then, given a regular neighborhood $N$ of $g^c$ in $M$, the 
collection of surfaces $\{ {\Sigma}_i\}$ is homotopic to $\{{\Sigma}'_i\}$ 
with homotopy supported in $N$, and $N$ contains an immersed disk ${\Delta}$ 
on ${\gamma}$ such that all surfaces $\Delta,{\Sigma}'_1,\ldots,
{\Sigma}'_{2^h-1}$ are pairwise disjoint. Moreover, the surfaces 
$\{ {\Sigma}_i\}$ stay pairwise disjoint during the homotopy.

The bound $2^h-1$ on the number of surfaces ${\Sigma}$, for which this 
conclusion holds in general, is sharp.
\end{thm}

The term {\em proper immersion} of a capped grope usually incorporates the 
condition that the body $g$ is embedded, and 
that the cap interiors are disjoint from $g$, so only cap--cap intersections
are allowed. This assumption is not needed in theorem \ref{surfaces}.
The necessary condition is that the surfaces $\Sigma$ may intersect only
the caps of $g^c$, but not the body $g$. 

Briefly, the idea of the proof is the following. Consider the special
case, when each body surface of $g^c$ has genus $1$ (thus $g^c$ has
$2^h$ caps), and when each cap intersects a single surface.
Then there must be two caps intersecting the same surface, say 
${\Sigma}_i$. We keep just these two caps for $g^c$, and using
surgery and contraction/pushoff, $g^c$ is made disjoint from
${\Sigma}_i$ as well. The general situation is reduced to this
special case via the operation of {\em grope splitting}, explained
in lemma \ref{refinement lemma}. 
As another application of grope splitting, we present a new proof of
the {\em Grope Lemma}:

\begin{thm} \label{grope lemma} \sl
Two $n$--component links in $S^3$ are link homotopic if and only if they cobound disjointly 
immersed annulus-like gropes of class $n$ in $S^3 \times I$.
\end{thm}

This result was originally stated in \cite{FL}, in the case when one of the links is trivial.
In the generality as stated here, Grope lemma is proved in \cite{KT}, using Milnor group. 
We refer the reader to \cite{FL}, \cite{KT} for the background, and for applications of this
result to the surgery conjecture. Note that the grope {\it class} corresponds to the index in
the lower central series of a group, while the grope {\it height} in theorems \ref{surfaces},
\ref{ecp} corresponds to the index in the derived series. Thus a grope of height $h$ has 
class $2^h$. To prove the Grope lemma, cap each grope by any transverse map of the disks
into $S^3\times I$. Note that a grope of class $n$ has $n$ caps, while there are
only $n-1$ gropes bounded by other link components. Now the argument using
grope splitting, identical to the proof of Theorem \ref{surfaces}, gives disjoint maps
of annuli (singular link concordance.) \endproof

The proof of theorem \ref{surfaces} also implies the result on ${\pi}_1$--null immersions:

\begin{thm}[Exponential contraction/pushoff 
{\cite[Theorem 3.5]{FT}}]\label{ecp}\nl Let ${\phi}\co{\pi}_1 G^c
\longrightarrow{\pi}$ be a group homomorphism with $(G^c,{\gamma})$ a
Capped Grope of height $h$.  If ${\phi}$ maps the double point loops
of $G^c$ to a set of cardinality at most $2^h-1$ in ${\pi}$ then $G^c$
contains a disk on ${\gamma}$ which is ${\pi}_1$--null under ${\phi}$.
\end{thm}

Note that while finding a ${\pi}_1$--null disk in theorem \ref{ecp} is similar
to finding a disk disjoint from other surfaces in theorem \ref{surfaces},
the converse -- showing that $2^h-1$ is the sharp bound in theorem~\ref{ecp} -- 
remains a central unsolved problem.
(If the bound $2^h-1$ in theorem~\ref{ecp} could be replaced by $2^h$, then 
one would find an embedded disk in the ``model'' capped gropes, cf \cite{FT}.)

In theorem \ref{ecp}, we allow cap--cap and cap--body intersections of $g^c$,
but it is important that the body $g$ on its own has no self-intersections. 
(More precisely, we allow only ${\pi}_1$--null self-intersections of the body.)
Also note that the proof goes through for any, not necessarily untwisted,
thickening $G^c$ of $g^c$.

We now make a brief digression to discuss the proof of Exponential
contraction/pushoff in \cite{FT}. For a given Capped Grope $G^c$ of
height $h$, the proof constructs, for each cap $C$, $2^h-1$ dual
spheres in $G^c$ with certain crucial disjointness properties. If
these dual spheres are used to resolve cap intersections, then the
grope is not left intact, but it has to be completely contracted,
using all its caps. More precisely, the dual spheres are built in the
``complete'' contraction. (Perhaps this point is not stated clearly in
the exposition of \cite{FT}.) It is the construction of these dual
spheres that requires developing the theory of colored link homotopy,
to show that a certain colored link $L$ is colored homotopically
trivial. In the present paper we do not follow that path, but prove
theorem \ref{ecp} directly. We note that our proof can also be used to
show that that particular link $L$ is trivial. (It also implies that
each dual sphere is {\em embedded}, so the intersections only occur
between {\em different} spheres of the same color.) The link $L$ is a
certain colored ramified iterated Bing double of the Hopf link, which
arises as a Kirby handle diagram of the Grope $G$. We say a few more
words about this link at the end of the proof of theorem
\ref{surfaces}. The theory of colored Milnor groups introduced in
\cite{FT} can be used under more general circumstances, for example
for the purpose of distinguishing non-homotopic links, but this
generality is not used in the proof of Exponential
contraction/pushoff.

\proof[Proof of Theorem \ref{surfaces}]
Consider first the special case when all body surfaces of $g$ have 
genus one, and each cap of $g^c$ intersects just one of the 
${\Sigma}_i$'s. This case captures the essence of the bound in 
theorems \ref{surfaces} and \ref{ecp}. Since there are $2^h$ caps and 
$2^h-1$ surfaces $\{ {\Sigma}_i\}$, at least two of the caps $C_1$, $C_2$ 
intersect the same surface ${\Sigma}_{i_0}$ (and they are disjoint from 
all other surfaces.) Consider these two caps, and for the rest of the proof
disregard all other caps of $g^c$. Suppose that ($C_1$, $C_2$) is a dual 
pair of caps, so they are attached to the symplectic pair of circles in 
an $h$-th stage surface of $g$. In this case contract $C_1$ and $C_2$ and 
push ${\Sigma}_{i_0}$ off the contraction to get ${\Sigma}'_{i_0}$.
Consider the disk ${\Delta}$ on ${\gamma}$ which ``uses'' only
the contraction of $C_1$ and $C_2$, and not the other caps. This
disk is gotten by successive surgeries along the branch of $g$ which leads 
from ${\gamma}$ to the tips $C_1$ and $C_2$; all other caps and surfaces 
in $g^c$ are disregarded. The disk ${\Delta}$ and the new surfaces
$\{{\Sigma}'_i\}$ satisfy the conclusion of the theorem. Note that
$g^c$ is not framed, so ``parallel copies'' (perturbations) of the surface
stages of $g^c$, which are used in the surgeries and contractions,
may intersect. This is not important in our argument, since the goal
is to find only an immersed disk.

If the caps $C_1$ and $C_2$ are not dual, still disregard all other caps,
surger two top stage surfaces, which are capped by $C_1$ and $C_2$
respectively, along these caps and continue surgering until the 
two new caps become dual. This reduces the situation to the previous case. 

Now consider the general case, with surfaces of an arbitrary genus, and 
when each cap may intersect several different
surfaces ${\Sigma}_i$. We will need for the proof the following operation
of {\em grope splitting}, so we make a digression to explain it in detail.

\begin{lemma}[Grope splitting] \label{refinement lemma} \sl
Let $(g^c,{\gamma})$ be a capped grope in $M^4$, and let
${\Sigma}_1,\ldots$, ${\Sigma}_n$ be surfaces in $M$, disjoint from the
body of $g^c$, but perhaps intersecting its caps. Then, given a regular neighborhood 
$N$ of $g^c$ in $M$, there is a capped grope $(g^c_{\rm split},{\gamma})
\subset N$, such that each cap of $g^c_{\rm split}$ intersects
at most one of the surfaces $\Sigma$, and each body surface, above
the first stage, of $g^c_{\rm split}$ has genus $1$.
\end{lemma}

\begin{figure}[ht]
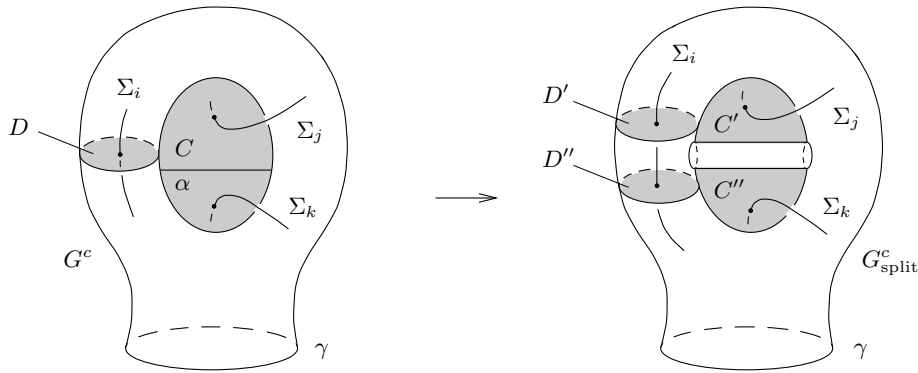

\small
\refine
\caption{Grope splitting}
\label{refinement}
\end{figure}

\proof First assume that $N$ is the untwisted thickening of $g^c$, $N=G^c$,
and moreover let $g^c$ be a model capped grope (without double points).
Let $C$, $D$ be a dual pair of its caps, and let ${\alpha}$ be an arc in 
$C$ with endpoints on the boundary of $C$. 
(In our applications, $\alpha$ will be chosen to separate intersection 
points of $C$ with different surfaces ${\Sigma}_j$, ${\Sigma}_k$, 
as shown in figure \ref{refinement}.) Recall that the untwisted 
thickening $N$ of $g^c$ is defined as the thickening in ${\mathbb{R}}^3$, 
times the interval $I$. We consider the $3$--dimensional thickening,
and surger the top-stage surface of $g$, which is capped by $C$ and $D$,
along the arc $\alpha$. The cap $C$ is divided by ${\alpha}$ into two 
disks $C'$, $C''$ which serve as the caps for the new grope; their dual 
caps $D'$, $D''$ are formed by parallel copies of $D$. This operation 
increases the genus of this top-stage surface by $1$; note that if some 
surface ${\Sigma}_i$ intersected the cap $D$ of $G^c$, it will intersect 
both caps $D'$, $D''$ of $G^c_{\rm split}$. 

We described this operation for a model capped grope; a splitting of a 
capped grope with double points is defined as the image of this operation 
in $N\!$/plum\-bings (where the arcs ${\alpha}$ are chosen to avoid the double 
points.) Also note that the same construction works for any (not
necessarily untwisted) thickening $N$: all that one needs is a line
subbundle of the normal bundle of the disk $C$ in $N$, restricted to
${\alpha}$. The fact that the new caps $D'$ and $D''$ may intersect
is not important here.

\begin{figure}[ht]
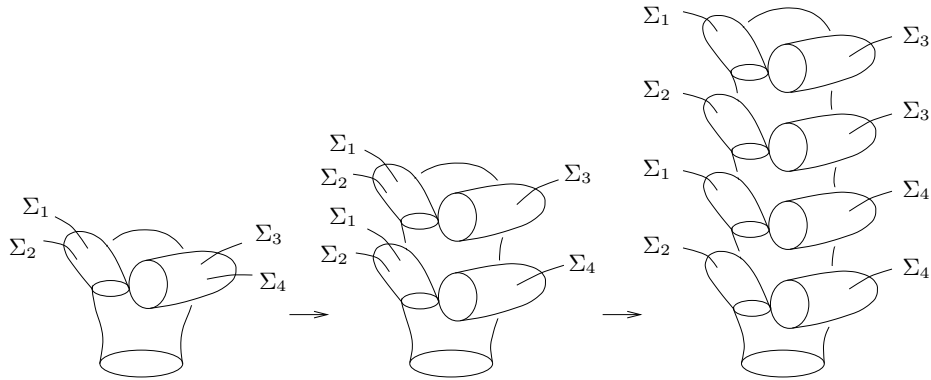

\small
\exrefine
\caption{Example of a splitting of a top stage surface of $g^c$}
\label{refinement example}
\end{figure}

Continue the proof of lemma \ref{refinement lemma} by dividing each cap 
$C$ by arcs $\{ {\alpha} \}$, so that each component of $C\smallsetminus
\cup{\alpha}$ intersects at most one surface in the collection $\{ 
{\Sigma}_i\}$, and splitting $g^c$ along all these arcs. (At most
$n$ arcs are needed for each cap.) The result is illustrated
in figure \ref{refinement example}. We apply the same operation to 
the surfaces in the $(h-1)$-st stage of the grope, separating each top 
stage surface by arcs into genus $1$ pieces. This procedure is performed
inductively, descending to the first stage of $g^c$. For example, if 
originally each body surface of $g$ had genus one, and each cap 
intersected all $n$ surfaces $\{ {\Sigma}_i\}$, then after this 
complete splitting procedure the first stage surface will have 
genus $n^{2^h}$. 
\endproof

\noindent
We continue the proof of theorem \ref{surfaces} in the general case,
applying this complete grope splitting procedure to $g^c$. Separate the 
first stage surface by arcs into genus one pieces and treat each one of 
them separately, as in the special (genus one) case, considered above. 
If one of the caps of $g^c$ is disjoint from all 
surfaces ${\Sigma}_i$, the result for that genus one piece follows trivially. 
The disk $\Delta$ bounding $\gamma$ is obtained as the union of disks 
produced by the genus one pieces.

\begin{figure}[ht]
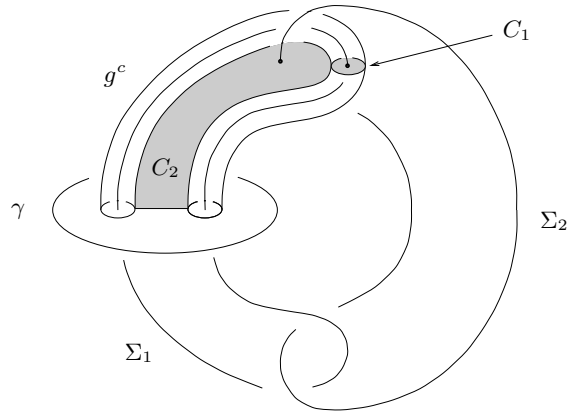

\small
\grope
\caption{$\gamma$ bounds a capped surface}
\label{capped grope}
\end{figure}

The proof that the bound $2^h-1$, for which the conclusion of the theorem 
holds in general, is sharp, is an elementary and well-known calculation 
in Massey products, or Milnor's $\bar\mu$--invariants. Consider the model 
capped grope $g^c$ (without double points) of height $h$ with each body 
surface of genus one and with the caps $C_1,\ldots,C_{2^h}$, and consider 
$2^h$ surfaces ${\Sigma}_1,\ldots,{\Sigma}_{2^h}$ such that for each $i$, 
the cap $C_i$ intersects ${\Sigma}_i$. The untwisted thickening of the 
model grope $g^c$ is the four--ball $B^4$; the attaching circle ${\gamma}$ 
of $g^c$ and the intersection of the surfaces ${\Sigma}_i$ with $\partial 
B^4=S^3$ form the Borromean rings, shown in Figure \ref{capped grope}, in 
the case $h=1$, cf \cite{Kirby}, \cite[12.2]{FQ}. The picture for larger $h$ is obtained by 
iterative Bing doubling (cf section 7 in \cite{K}) of the link components, 
other than $\gamma$. At each step of the iteration, the caps are replaced 
by genus one capped surfaces (copies of figure \ref{capped grope}.) 
The components of the resulting link do not bound disjoint maps 
of disks in $B^4$ since any iterated Bing double of the Hopf link 
has non-trivial $\bar\mu$--invariants \cite{Milnor}.
\endproof

\proof[Proof of Theorem \ref{ecp}] This is similar to the proof of
Theorem \ref{surfaces} above, only instead of separating intersections
of the caps with different surfaces, the grope splitting procedure
will now be used to separate intersection points among the caps of
$g^c$, which correspond to different group elements in $\pi$. Recall
from \cite[2.9]{FQ} that the new group element created by the
operation of pushoff from the elements $f$ and $g$ is $f\cdot g^{-1}$,
thus only trivial double point loops are created during the final
step.

When one applies the grope splitting in order to separate the
selfintersections of $g^c$, rather than intersections of $g^c$ with other
surfaces, one cannot achieve the situation shown in figure \ref{refinement
example} where each cap has precisely one double point. Indeed, splitting a 
cap $C$ requires using two copies of the dual cap $D$, 
making it impossible to achieve progress in this respect. However, the double 
point loops, produced by the parallel copies $D'$, $D''$ of $D$, give the same 
group element, and the new intersections $D'\cap{\Sigma}_i$, $D''\cap{\Sigma}_i$ 
do not need to be separated in ${\Sigma}_i$.  

Another subtlety concerns ordering the sheets at each intersection point
(if one ordering gives an element $g$ in $\pi$, switching them
gives $g^{-1}$.) The statement of theorem \ref{ecp} implicitly contains
a choice of the first sheet at each intersection point. The proof
of theorem \ref{surfaces}, followed here without a change, would
only give the bound $2^{h-1}-1$. Thus a slight correction is necessary in 
the situation after the grope splitting is completed, when two caps 
$C_1$ and $C_2$ on the same branch $B$ have intersections
with some other caps, representing the same non-trivial element $g$, but where
$C_1$ is considered as the first sheet for its intersection point, $C_2$ 
is considered as the second sheet for its intersection, and $g\neq g^{-1}$. 
If both elements $g$, $g^{-1}$ 
are on the list of $2^h-1$ elements, then this problem does not arise.
Suppose $g$ is on the list, and $g^{-1}$ is not. Take the cap $C_1$
(labeled as the first sheet) and surger the grope along the branch leading 
to $C_1$ (all other caps, including $C_2$, of this branch are discarded.) Let $C$ be a 
cap, lying on some branch $B'$, intersecting $C_1$, giving the group 
element $g$. Note that $C$ is considered as the second sheet for this intersection, 
and it will be discarded when this operation is applied to the corresponding branch $B'$
of $g^c_{\rm split}$. Hence this procedure is consistent, giving raise at the end to a ${\pi}_1$--null 
disk~$\Delta$.

Note that theorem \ref{ecp} holds for any (not necessarily untwisted)
four--dimensional thickening of $g^c$. The ``parallel copies'' of the
surfaces, that have to be taken for surgeries and contractions in the 
proof, are just perturbations of the originals. The resulting 
singularities are acceptable, since their double point loops are 
trivial in ${\pi}_1$.
\endproof

{\bf Acknowledgements}\qua Research was supported by the Institute for
Advanced Study (NSF grant DMS 97-29992) and by Harvard University.

\end{document}